
\count100= 1
\count101= 10

\overfullrule 0pt

\magnification\magstep1
\hsize 4.43in
\vsize 7.28in



\font\bfab = cmbx9
\font\fa=cmr17
\font\fb=cmr12

\font\fab=cmr9
\font\sfab=cmr7
\font\fs=cmr6
\font\fd=cmr5
\font\slab=cmsl9
\font\tab= cmmi9
\font\sab= cmmi7
\font\ssab= cmmi6
\font\abst= cmsy9
\font\sabs= cmsy7
\font\ssabs= cmsy6
\font\itab= cmti9

\def\figfont{
    \textfont0 = \fab
    \scriptfont0 = \sfab
    \scriptscriptfont0 = \fs
    \textfont1 = \tab
    \scriptfont1 = \sab
    \scriptscriptfont1 = \ssab
    \textfont2 = \abst
    \scriptfont2 = \sabs
    \scriptscriptfont2 = \ssabs
    \let \sl = \slab
    \let \bf = \bfab
    \let \it = \itab
    \baselineskip 9pt
    \fab}


\def\CC{{\rm C\kern-.18cm\vrule width.6pt height 6pt depth-.2pt
\kern.18cm}}

\def\NN{{\mathop{{\rm I}\kern-.2em{\rm N}}\nolimits}}

\def\PP{{\mathop{{\rm I}\kern-.2em{\rm P}}\nolimits}}

\def\RR{{\mathop{{\rm I}\kern-.2em{\rm R}}\nolimits}}

\def\RRt{{\fa I}\kern-.2em{\fa R}}


\def\ZZ{{\mathop{{\rm Z}\kern-.28em{\rm Z}}\nolimits}}





\def\makebold#1{\mathord{\setbox0=\hbox{#1}%
       \copy0\kern-\wd0%
       \raise\dimen1\copy0\kern-\wd0%
       {\advance\dimen1 by \dimen1\raise\dimen1\copy0}\kern-\wd0%
       \kern\dimen0\raise\dimen1\copy0\kern-\wd0%
       {\advance\dimen1 by \dimen1\raise\dimen1\copy0}\kern-\wd0%
       \kern\dimen0\raise\dimen1\copy0\kern-\wd0%
       {\advance\dimen1 by \dimen1\raise\dimen1\copy0}\kern-\wd0%
       \kern\dimen0\raise\dimen1\copy0\kern-\wd0%
       \kern\dimen0\box0}}


\newbox\maboite

\def\boxit#1#2{\setbox\maboite=\hbox{\kern#1{#2}\kern#1}%
    \dimen1=\ht\maboite \advance\dimen1 by #1 \dimen2=\dp\maboite
\advance\dimen2 by #1%
    \setbox\maboite=\hbox{\vrule height\dimen1
depth\dimen2\box\maboite\vrule}%
    \setbox\maboite=\vbox{\hrule\box\maboite\hrule}%
    \advance\dimen1 by .4truept \ht\maboite=\dimen1%
    \advance\dimen2 by .4truept \dp\maboite=\dimen2 \box\maboite\relax}




%

\def\frac#1#2{{#1 \over #2}}


\def\ms{\medskip}

\def\noin{\noindent}










\def\pf{\noindent{\bf Proof: }}

\def\eop{\makeblanksquare6{.4}\ms}

\def\makeblanksquare#1#2{
\dimen0=#1pt\advance\dimen0 by -#2pt
      \vrule height#1pt width#2pt depth0pt\kern-#2pt
      \vrule height#1pt width#1pt depth-\dimen0 \kern-#1pt
      \vrule height#2pt width#1pt depth0pt \kern-#2pt
      \vrule height#1pt width#2pt depth0pt
}



\def\abstract#1{\bigskip\bigskip\medskip%
    {\narrower \baselineskip 9pt \fab \noindent {\bf Abstract.~~}%
    \textfont0 = \fab
    \scriptfont0 = \sfab
    \scriptscriptfont0 = \fs
    \textfont1 = \tab
    \scriptfont1 = \sab
    \scriptscriptfont1 = \ssab
    \textfont2 = \abst
    \scriptfont2 = \sabs
    \scriptscriptfont2 = \ssabs
    \let \it = \itab
    \let \sl = \slab
    #1\bigskip}\medskip}
\def\author#1{\bigskip\bigskip\centerline{\fb #1}}

\def\copyright{\hbox{{\fb o}\kern-.61em \raise .46ex \hbox{\fd c}}}
\def\footnoterule{\kern -3pt \hrule width 0truein \kern 2.6pt}
\def\leftheadline{\ifnum\pageno=\count100 \hfill%
  \else\rm\folio\hfil\it\shortauthor\fi}
\def\rightheadline{\ifnum\pageno=\count100 \hfill%
  \else\it\shorttitle\hfil\rm\folio\fi}
\def\title#1{\centerline {\fa #1}}

\def\titexp#1#2{\hbox{{\fa #1} \kern-.25em \raise .90ex \hbox{\fb #2}}\/}
\def\titsub#1#2{\hbox{{\fa #1} \kern-.25em \lower .60ex \hbox{\fb #2}}\/}

\nopagenumbers
\headline{\ifodd\pageno\rightheadline \else\leftheadline\fi}
\footline{\hfil}
\null\vskip 18pt
\centerline{}
\pageno=\count100
\count102=\count100
\advance\count102 by -1
\advance\count102 by \count101




\def\sect#1{\goodbreak\bigskip\smallskip\centerline{\bf\S #1}\medskip
    \noindent\ignorespaces}


\def\subsect#1{\goodbreak\bigskip\leftline{\bf#1}\medskip
  \noindent\ignorespaces}




\def\Address{\nonfrenchspacing\goodbreak\bigskip\obeylines}

\def\Remark#1.{\goodbreak\medskip\noin {\bf Remark#1.}}
\def\Example#1.{\goodbreak\medskip\noin {\bf Example#1.}}


\def\ref{\smallskip\global\advance\refnum by 1 \item{\the\refnum.}}
\newcount\refnum \refnum = 0

\def\References{\goodbreak\bigskip\centerline{\bf References}\bigskip
   \frenchspacing}


\def\CA{Constr.\ Approx.}

\def\CRASP{C. R. Acad.\ Sci.\ Paris}

\def\NA{Numer.\ Algorithms}



\title{Generalized $C^1$ quadratic B-splines}
\title{generated by Merrien subdivision algorithm}
\title{and some applications}

\author{ Paul Sablonni\`ere}

\def\shorttitle{Generalized quadratic B-splines}

\def\shortauthor{ps}

\abstract{
A new global basis of B-splines is defined in the space of generalized
quadratic splines (GQS)
generated by Merrien subdivision algorithm. Then, refinement equations  for
these B-splines and
the  associated corner-cutting algorithm are given. Afterwards, several
applications are
presented. First a global  construction of monotonic and/or convex
generalized splines interpolating  monotonic and/or convex data. Second,
convergence of sequences
of control polygons to the graph of a GQS. Finally, a Lagrange  interpolant
and a quasi-interpolant which are exact on the space of affine polynomials
and  whose infinite
norms are uniformly bounded independently of the partition.}

*********************************

\sect{1. Introduction and notations}

\noin
This paper is a continuation of [7] and [10] where a one-parameter family
of $C^1$ Hermite
interpolants was defined by the Merrien subdivision algorithm ([6], abbr. MSA).
This family provides a general solution to arbitrary monotone and convex
Hermite interpolation
problems with data at the ends of a single interval. Here we extend and
complete these results by
considering Hermite interpolation on a partition of a given interval by
$C^1$ functions whose
restrictions to subintervals are generated by MSA.
Let $X:=\{a=x_0,x_1,\ldots,x_n=b\}$ be such a partition of $I=[a,b]$ with
$I_i=[x_{i-1},x_i]$ and $h_i=x_i-x_{i-1}$ for $1\le i\le n$.
Let ${\beta}=(\beta_1,\ldots\beta_n)$ be a sequence of parameters
$\beta_i\in [-1,0[$.

\noin
On each subinterval $I_i$, we consider the $4$-dimensional space $V(\beta_i)$
generated by Merrien subdivision algorithm and depending on the specific
pair of parameters
$(\alpha_i,\beta_i)$, where
$\alpha_i=\frac{\beta_i}{4(1-\beta_i)}\in [-\frac18,0[$.
In that case, it was proved that the MSA is
$C^1$-convergent. Moreover, for $f\in V(\beta_i)$, $f'$ is H\"older and
satisfies the following
inequality for some constant $C>0$ and  for
$\gamma_i=-log_2(1+\frac12\beta_i)$ (see [6], proposition 2):
$$
\vert f'(x)-f'(y) \vert \le C\vert x-y\vert^{\gamma_i},\;\; (x,y)\in
I_i\times I_i.
$$
We denote by $GS_2(I,X,\beta)$ or simply $GS_2(\beta)$ the subspace of all
functions
$g\in C^1(I)$ whose restriction to
$I_i$ is in $V(\beta_i)$ for $1\le i\le n$. The dimension of this space is
$2n+2$ since $g$
is uniquely determined by Hermite data $y_i=g(x_i)$ and $p_i=g'(x_i)$ for
$0\le i\le n$.
The elements of $GS_2(\beta)$ are called {\sl generalized quadratic
splines} (abbr.
GQS) because their  properties are similar to those of classical quadratic
splines, which
correspond to the choice $\beta_i=-1, 1\le i\le n$.

\noin
Here is an outline of the paper. In section 2, we define a {\sl new global
basis of  B-splines}
for the space $GS_2(\beta)$ of GQS. In section 3, we use these B-splines to
express in a new
global form the results on {\sl monotonicity and convexity preserving
properties} of GQS, already
given in [7][10] for functions defined on one subinterval. In section 4, we
also give a global
version of the algorithms constructing monotonic or convex interpolants
given in the same papers.
When compared to other $C^1$ Hermite interpolants of the literature, the
advantage of our
method is its simplicity both in the local construction of the GQS and in
the adaptivity of the
method to arbitrary sets of data. For example, in [1][2], the method is
adaptive, however it may
need polynomials of arbitrary high degree,
which can lead to rather complicated expansions and calculations.

\noin
In section 5, we give the {\sl refinement equation} for coarse B-splines in
terms of fine
B-splines in the space of GQS defined on a refinement of the initial partition. From this result
we deduce a {\sl corner-cutting algorithm} (see [9] for definitions and
properties) which is the
geometric form of the algorithm expressing the
new coefficients of a GQS $g$ in the B-spline basis of the fine space in
terms of its old
coefficients in the B-spline basis of the coarse space. We also prove the
convergence to the
graph of $g$ of the sequence of control polygons associated with successive
steps of this
algorithm.

\noin
Finally, in section 6, we study two approximation operators: a {\sl
Lagrange interpolant} and a
{\sl quasi-interpolant} which are both exact on the space $\PP_1$ of affine
polynomials and whose
uniform norm are {\sl uniformly bounded independently of the given
partition} on $I$. This
extends a previous result given for ordinary quadratic splines by Kammerer,
Reddien
and Varga [4] and also by Marsden [5]. We postpone numerical applications
to a further paper
which should contain variants and refined versions of the general
algorithms presented in
Section 4 of the present paper.

\ms
\noin
For the sake of clarity, we now recall the basic equations of the MSA
giving the values at the
midpoint $m=\frac12(a+b)$ of $[a,b]$ of a function $f$ and its first
derivative  $f'$
from the four values:
$$
\{f(a),f'(a);f(b),f'(b)\}
$$
at the endpoints of the interval (see [5]).
The construction starts with $[a,b]=I_i=[x_{i-1},x_i]$  and gives the values of
$f$ and $f'$ at the dyadic points of $I_i$. Let $h=b-a$ and
$\theta_i=\frac12\frac{\beta_i}{\beta_i-1}=-2\alpha_i\in ]0,\frac14]$, then
$$
f(m)=\frac12\left((f(a)+f(b))-\theta_i h (f'(b)-f'(a))\right)
$$
$$
f'(m)=\frac{1}{1-2\theta_i}\left(\frac{f(b)-f(a)}{h}-2\theta_i
\frac{f'(b)+f'(a)}{2}\right)
$$

\noin
In each  subinterval $I_i$, let us define the two points
$$
\xi_i=x_i-\theta_ih_i, \quad \eta_i=x_i+\theta_{i+1}h_{i+1},
$$
with $\xi_0=x_0$ and $\eta_n=x_n$.
Then each element $g_i\in V(\beta_i)$ can be expressed as
$$
g_i=a_{i-1}b_0+d_{i-1}b_1+c_ib_3+a_ib_4
$$
in the {\sl local B-spline basis} $\{b_0,b_1,b_2,b_3\}$ of $V(\beta_i)$
defined in [10].
By definition, the quadruplet  $[a_{i-1}, d_{i-1}, c_i, a_i]$ is the list
of {\sl
B-coefficients} of
$g_i$ on the subinterval $I_i$. The {\sl local control polygon} (abbr. LCP)
of $g_i$ has
the four following {\sl local control vertices}
$$
\tilde a_{i-1}=(x_{i-1},a_{i-1}), \quad \tilde
d_{i-1}=(\eta_{i-1},d_{i-1}),\quad
\tilde c_i=(\xi_i,c_i), \quad \tilde a_i=(x_i,a_i).
$$
The ordinates of these points are the B-coefficients of $g_i$ and they can
be expressed in
function of the four {\sl Hermite data} $(y_{i-1},p_{i-1};y_i,p_i)$ at the
two end-points of
$I_i$:
$$
y_i=g_i(x_i)=g_{i+1}(x_i), \quad p_i=g'_i(x_i)=g'_{i+1}(x_i).
$$
They are given, for $0\le i \le n-1$, by
$$
a_i=y_i, \quad c_i=a_i-\theta_ih_ip_i, \quad
d_i=a_i+\theta_{i+1}h_{i+1}p_i, \quad a_{i+1}=y_{i+1}
$$
(with the convention $h_0=h_{n+1}=0$).

\noin
Using the MSA and the properties
$\alpha_i=-\frac12\theta_i$,
$\beta_i=\frac{2\theta_i}{2\theta_i-1}$, $1-\beta_i=\frac{1}{1-2\theta_i}$, and
$\xi_i-\eta_{i-1}=(1-2\theta_i)h_i$, we obtain respectively at the midpoint
$\displaystyle{m_i=\frac{x_{i-1}+x_i}{2}}$ of $I_i$:
$$
\displaylines{
g_i(m_i)=\frac12\left((a_{i-1}+a_i)-\theta_i h_i(p_i-p_{i-1})\right)
=\frac12(d_{i-1}+c_i),\cr
g'_i(m_i)=\frac{1}{1-2\theta_i}\left(\frac{a_i-a_{i-1}}{h_i}-2\theta_i\frac{p_{i-1}+p_i}{2}\right)
=\frac{c_i-d_{i-1}}{\xi_i-\eta_{i-1}},
}
$$
which proves that the tangent to the curve at $m_i$ is the segment $\tilde
d_{i-1}\tilde c_i.$

*********************************

\sect{2. Generalized quadratic B-splines}

\noin
In this section, we define a {\sl global B-spline basis} of the space
$GS_2(\beta)$
of generalized quadratic splines.

\noin
Let
$\displaystyle{\omega_i=\frac{\theta_{i+1}h_{i+1}}{\theta_ih_i+\theta_{i+1}h
_{i+1}}}$
for $1\le i\le n-1$.

\ms
\noin
{\bf Definitions}

\noin
(i) For $1\le i\le n-1$, let $B_{2i}$ be the function whose support
is the segment $[x_{i-1},m_{i+1}]$  and whose lists of B-coefficients in
the subintervals
$I_i$ and $I_{i+1}$ are respectively
$$
[0,0,1,\omega_i], \quad [\omega_i,0,0,0].
$$
(ii) Similarly, let $B_{2i+1}$ be the function whose support
is the segment $[m_i,x_{i+1}]$, and whose lists of B-coefficients in the
subintervals
$I_i$ and $I_{i+1}$ are respectively
$$
[0,0,0,1-\omega_i], \quad [1-\omega_i,1,0,0].
$$
(iii) Moreover, there are $4$ special B-splines
$\{B_0,B_1,B_{2n},B_{2n+1}\}$ at the end points
of $I$, defined respectively by

\noin
supp($B_0$)=$[x_0,m_0]$, its list of B-coefficients in $I_1$ being $[1,0,0,0]$.

\noin
supp($B_1$)=$[x_0,x_1]$, its list of B-coefficients in $I_1$ being $[0,1,0,0]$.

\noin
supp($B_{2n}$)=$[x_{n-1},x_n]$, its list of B-coefficients in $I_n$ being
$[0,0,1,0]$.

\noin
supp($B_{2n+1}$)=$[m_n,x_n]$, its list of B-coefficients in $I_n$ being
$[0,0,0,1]$.
\ms
\noin
It is easy to verify that
$$
B'_0(x_0)=\frac{-1}{\eta_0-x_0}=\frac{-1}{\theta_1 h_1}, \quad
B'_{2n}(x_n)=\frac{-1}{x_n-\xi_n}=\frac{-1}{\theta_n h_n},
$$
$$
B'_1(x_0)=-B'_0(x_0), \quad B'_{2n+1}(x_n)=-B'_{2n}(x_n),
$$
and for $1\le i\le n-1$
$$
B'_{2i}(x_i)=\frac{-1}{\eta_i-\xi_i}=-\frac{\omega_i}{\theta_{i+1}h_{i+1}}=-
B'_{2i+1}(x_i).
$$

\proclaim Theorem 1. The generalized quadratic B-splines $\{B_k, 0\le k\le
2n+1\}$ form a basis of
the space $GS_2({\beta})$. Moreover, they form a partition of unity, or a
blending system, in
this space.

\pf
First, let us prove that the B-splines belong to the space $GS_2({\beta})$,
i.e. that they
are  $C^1$ continuous at the points $\{x_1,x_2,\ldots,x_{n-1}\}$. By
construction, they are
already $C^1$ in each subinterval $I_i$. In addition, their derivatives
satisfy the above
relations, so they are continuous at the interior points of $X$.

\noin
The B-splines are linearly independent: assume that $g=\sum_{k=0}^{2n+1}
\gamma_k B_k=0$.
Then, the B-coefficients of the restriction $g_i=\sum_{k=2i-2}^{2i+1}
\gamma_k B_k$ of $g$ to
$I_i$  are respectively
$$
a_{i-1}=\omega_{i-1}\gamma_{2i-2}+(1-\omega_{i-1})\gamma_{2i-1}, \quad
d_{i-1}=\gamma_{2i-1},
$$
$$
c_i=\gamma_{2i}, \quad a_i=\omega_i\gamma_{2i}+(1-\omega_i)\gamma_{2i+1}.
$$
Therefore, since $\omega_{i-1}$ and $1-\omega_i$ are non zero, we get
successively
$\gamma_{2i-1}=\gamma_{2i}=0$, and $\gamma_{2i-2}=\gamma_{2i+1}=0$.

\noin
The  B-splines generate the space $GS_2({\bf\beta})$: it suffices to express the
coefficients $\gamma_k$ in function of the B-coefficients.
The above equations give immediately
$$
\gamma_{2i}=c_i, \quad \gamma_{2i+1}=d_i.
$$
Finally, let us prove that $\sum  B_k=1$: it suffices to prove it on each
subinterval $I_i$. As
local  B-spline bases are blending systems, one has to show that the sum of
local B-coefficients
of global B-splines is equal to $1$. This property is easily deduced from
the lists
of B-coefficients of
$B_{2i-2},B_{2i-1},B_{2i},B_{2i+1}$ on the interval $I_i$ which are
respectively equal to
$$
[\omega_{i-1},0,0,0],\;\;[1-\omega_{i-1},1,0,0],\;\;[0,0,1,\omega_i],\;\;[0,
0,0,\omega_i]
$$
\eop

*********************************

\sect{3. Global and local control polygons. Monotonicity and convexity}

\noin
It has been proved in [6] that the spaces $V(\beta_i)$ contain the space
$\PP_1$ of
affine polynomials . As the list of local B-coefficients on $I_i$  of the
function $e_1(x)=x$
is $[x_{i-1}, \eta_{i-1},\xi_i,x_i]$, the results of section 2 show that the
representation of $e_1$ in the basis of generalized B-splines is given by
$$
e_1=\sum_{i=0}^{n+1} (\xi_i B_{2i}+\eta_i B_{2i+1}).
$$
Therefore, we can define the global {\sl spline control polygon} (abbr. SCP) of
$g=\sum_{k=0}^{2n+1}\gamma_k B_k$ with vertices
$$
\tilde \gamma_{2i}=(\xi_i, \gamma_{2i}), \quad \tilde
\gamma_{2i+1}=(\eta_i, \gamma_{2i+1})
$$
from which we easily deduce the vertices of the local control polygon (LCP)
$$
\tilde a_{i-1}=\omega_{i-1}\tilde
\gamma_{2i-2}+(1-\omega_{i-1})\tilde\gamma_{2i-1}, \quad
\tilde c_{i-1}=\tilde \gamma_{2i-2},
$$
$$
\quad \tilde d_i=\gamma_{2i+1},
\quad \tilde a_{i}=\omega_{i}\tilde
\gamma_{2i}+(1-\omega_{i})\tilde\gamma_{2i+1}.
$$
\proclaim Theorem 2.
A function $g\in GS_2(\beta)$ is monotonic (resp. convex) if and only if
its SCP is
monotonic (resp. convex), with the same sense of variation.

\pf
>From theorem 6 of [10], we know that $g$ is (e.g.) increasing on $I_i$ if
and only if
its LCP is increasing, i.e. iff
$$
a_{i-1}\le d_{i-1}\le c_i\le a_i.
$$
Since we have
$$
d_{i-1}-a_{i-1}=\omega_{i-1}(\gamma_{2i-1}-\gamma_{2i-2)}, \;\;\;
c_i-d_{i-1}=\gamma_{2i}-\gamma_{2i-1},
$$
$$
a_i-c_i=(1-\omega_i)(\gamma_{2i+1}-\gamma_{2i}),
$$
we see that the above inequalities are satisfied iff $\gamma_{2i-1}\le
\gamma_{2i-1}\le \gamma_{2i-1}$, i.e.
iff the global SCP is increasing.

\noin
Similarly, from theorem 8 of [10], we know that $g$ is convex on $I_i$ if
and only if
$$
\frac{d_{i-1}-a_{i-1}}{\eta_{i-1}-x_{i-1}}\le\frac{c_i-d_{i-1}}{\xi_i-\eta_{
i-1}}\le \frac{a_i-c_i}{x_i-\xi_i}.
$$
Using the equalities
$$
\eta_{i-1}-\xi_{i-1}=\frac{\theta_ih_i}{\omega_{i-1}},\;\;
\xi_i-\eta_{i-1}=(1-2\theta_i)h_i,\;\;
\eta_i-\xi_i=\frac{\theta_ih_i}{1-\omega_i}
$$
we obtain successively the following identities
$$
\frac{d_{i-1}-a_{i-1}}{\eta_{i-1}-x_{i-1}}=\frac{\omega_{i-1}(\gamma_{2i-1}-
\gamma_{2i-2})}
{\theta_ih_i}=\frac{\gamma_{2i-1}-\gamma_{2i-2}}{\eta_{i-1}-\xi_{i-1}},
$$
$$
\frac{c_i-d_{i-1}}{\xi_i-\eta_{i-1}}=\frac{\gamma_{2i}-\gamma_{2i-1}}{(1-2\theta_i)h_i}=
\frac{\gamma_{2i}-\gamma_{2i-1}}{\xi_i-\eta_{i-1}},
$$
$$
\frac{a_i-c_i}{x_i-\xi_i}=\frac{(1-\omega_i)(\gamma_{2i+1}-\gamma_{2i})}{\theta_ih_i}=
\frac{\gamma_{2i+1}-\gamma_{2i}}{\eta_i-\xi_i},
$$
which show that the convexity of the global SCP is equivalent to that of
the LCP in each
subinterval, i.e. to the convexity of $g$.
\eop

*************************

\sect{4. Construction of monotone or/and convex interpolants}

\subsect{Algorithm 1:  monotone interpolant}

\noin
This algorithm describes the construction of an increasing interpolant

\noin
$g\in GS_2(I,X,\beta)$ to arbitrary increasing Hermite data:
$$
(y_i,p_i),\;\; 0\le i\le n,
$$
assumed to satisfy the properties:
$$
\Delta y_i > 0, \;\; p_i> 0.
$$
(There is a similar algorithm for decreasing data).
The main point consists in choosing the sequence $\beta$ in function of the
data.
>From theorem 2, we know that $g=\sum_k \gamma_k B_k$ is increasing if and
only if the sequence
$(\gamma_k)$ of its S-coefficients is increasing. On one hand, we already
know that
$$
\gamma_{2i+1}-\gamma_{2i}=d_i-c_i=(\theta_{i+1}h_{i+1}+\theta_i h_i)p_i> 0.
$$
On the other hand, using the notations:
$$
\tau_i=\frac{\Delta y_{i-1}}{h_i},\;\; \mu_i=\frac12 (p_{i-1}+p_i),
$$
we must have
$$
\gamma_{2i}-\gamma_{2i-1}=c_i-d_{i-1}=h_i(\tau_i-2\theta_i \mu_i)>0,
$$
which implies the following condition on the parameter
$\theta_i=\frac12{\beta_i}{\beta_i-1}:$
$$
0<\theta_i< \bar\theta_i=\frac12 \frac{\tau_i}{\mu_i}.
$$
Note that this condition is  equivalent to $g'(m_i)>0$.

\ms
\noin
There appear two cases and we obtain the following

\proclaim Theorem 3.
(i) if $\mu_i\le 2\tau_i$, then we can choose $\theta_i=\frac14$: in that
case, the local interpolant
in the subinterval $I_i$ is an ordinary quadratic spline ($\beta_i=-1$).
\ms
(ii) if  $\mu_i> 2\tau_i$, then we have to choose
$\theta_i\le\bar\theta_i=\frac12 \frac{\tau_i}{\mu_i}<\frac14$:
in that case, the local interpolant in the subinterval $I_i$ is a
generalized quadratic spline associated
with the parameter $\beta_i=\frac{-2\theta_i}{1-2\theta_i}$.

\subsect{Algorithm 2: convex interpolant}

\noin
This algorithm describes the construction of a convex interpolant

\noin
$g\in GS_2(I,X,\beta)$ to arbitrary convex Hermite data:
$$
(y_i,p_i),\;\; 0\le i\le n,
$$
assumed to satisfy the properties
$$
p_{i-1}<\tau_i<p_i
$$
on each subinterval $I_i$.
(There is a similar algorithm for concave data).
\noin
Using the results given in the proof of theorem 2, we get the convexity
conditions on $I_i$ for
the parameter $\theta_i=\frac12{\beta_i}{\beta_i-1}:$
$$
p_{i-1}\le \frac{\tau_i-2\theta_i \mu_i}{(1-2\theta_i)}\le p_i
$$
which can also be written
$$
(1-\theta_i)p_{i-1}+\theta_i p_i\le \tau_i\le (1-\theta_i)p_{i}+\theta_i p_{i-1}
$$
\noin
As for algorithm 1, there appear two cases:

\proclaim Theorem 4.

(i) if $\tau_i$ satisfies the two inequalities

$\frac14(3 p_{i-1}+p_i)\le \tau_i\le \frac14 (p_{i-1}+3 p_i)$,
then we chose $\theta_i=\frac14$: in that case, the local interpolant is an
ordinary quadratic spline.
\ms
(ii) else the local interpolant is a generalized quadratic spline
($\theta_i<\frac14$).
\ms
a) either $p_{i-1}<\tau_i<\frac14(3 p_{i-1}+p_i)$: then we have to choose
$\theta_i\le\bar\theta_i=\frac{\tau_i-p_{i-1}}{p_i-p_{i-1}}$.
\ms
b) or $\frac14 (p_{i-1}+3 p_i)<\tau_i<p_i$: then we have to choose
$\theta_i\le\bar\theta_i=\frac{p_i-\tau_i}{p_i-p_{i-1}}$.

In both cases, the GQS belongs to $V(\beta_i)$, with
$\beta_i=\frac{-2\theta_i}{1-2\theta_i}$.

\subsect{Algorithm 3: monotone and convex interpolant}

\noin
For sake of simplicity, we assume that the data are increasing and convex
$$
0<p_{i-1}<\tau_i<p_i
$$
on each subinterval $I_i$.
(There are  similar algorithms for the three other cases).
By putting together conditions of theorems 3 and 4, we get the same
algorithm as in the preceding case.

\proclaim Theorem 5.
(i) \quad if $\frac14(3 p_{i-1}+p_i)\le \tau_i\le \frac14 (p_{i-1}+3 p_i)$,
then we can choose $\theta_i=\frac14$.
The local interpolant is an ordinary quadratic spline.
\ms
(ii)  \quad else the local interpolant is a generalized quadratic spline.
\ms
a) either $\tau_i<\frac14(3 p_{i-1}+p_i)$, then we chose
$\theta_i\le\bar\theta_i=\frac{\tau_i-p_{i-1}}{p_i-p_{i-1}}$.
\ms
b) or $\tau_i>\frac14(p_{i-1}+3 p_i)$, then we chose
$\theta_i\le\bar\theta_i=\frac{p_i-\tau_i}{p_i-p_{i-1}}$.

In both cases, the GQS belongs to $V(\beta_i)$, with
$\beta_i=\frac{-2\theta_i}{1-2\theta_i}$.

\noin
Remark: these are only rough algorithms: in practice, one has to smooth a
little bit the above conditions
and also to treat the cases when there appear equalities in the conditions.
This will be done in a further more complete paper illustrated with
numerical examples.

algorithm ************************************

\sect{5. Refinement  equations and global corner-cutting algorithm}

\noin
In this section, we consider the subpartition $\bar{X}=X\cup\{m_i, 1\le
i\le n\}$ dividing $I$
into $2n$ subintervals. The space $GS_2(I, X,\beta)$ is a subspace of
dimension $2n+2$ of the new
space
$GS_2(I,\bar X,\bar \beta)$ of dimension $4n+2$. Here $\bar \beta$ denotes
the sequence of
parameters deduced from $\beta$ by taking twice the same parameter
$\beta_i$, once for the left
subinterval $I'_i=[x_{i-1},m_i]$ and once for the right subinterval
$I''_i=[m_i,x_i]$ of $I_i$:
$$
\bar \beta=(\beta_1, \beta_1,\beta_2,\beta_2,\ldots,\beta_n,\beta_n).
$$
The finer B-splines of $GS_2(I,\bar X,\beta)$ are denoted
$$
\{\bar B_l, \; 0\le l\le 4n+1\}.
$$
Thanks to the local corner-cutting algorithm (abbr. CCA, see [8], section
), one can compute the local
B-coefficients in the subintervals $I'_i$ and $I''_i$ of  the coarser
B-splines $\{B_k, 0\le k\le n+1\}$
in function of  their previous B-coefficients in $I_i$, and
we get the following refinement equations.

\noin
Let us recall that
$\displaystyle{\omega_i=\frac{\theta_{i+1}h_{i+1}}{\theta_ih_i+\theta_{i+1}h
_{i+1}}}$, for $1\le i\le n-1$,
$\omega_0=1-\omega_n=1.$

\proclaim Theorem 6.
(i) For all $1\le i\le n-1$, one has the two following refinement equations:
$$
B_{2i}=(\frac12+\frac14 \beta_i) \bar B_{4i-2}+(\frac12-\frac14 \beta_i)
\bar B_{4i-1}
+\frac12 (1+\omega_i) \bar B_{4i}+\frac12 \omega_i \bar B_{4i+1},
$$
$$
B_{2i+1}=\frac12 (1-\omega_i) \bar B_{4i}+\frac12 (2-\omega_i) \bar B_{4i+1}
+(\frac12-\frac14 \beta_i) \bar B_{4i+2}+(\frac12+\frac14 \beta_i) \bar
B_{4i+3}.
$$
(ii) For the B--splines at the endpoints, one has respectively:
$$
B_0=\bar B_0+\frac12\bar B_1, \;\; B_1=\frac12\bar B_1+(\frac12-\frac14\beta_1) \bar B_{2}+(\frac12+\frac14 \beta_1) \bar B_{3},
$$
$$
B_{2n}=\frac12\bar B_{4n}+(\frac12-\frac14 \beta_n) \bar
B_{4n-1}+(\frac12+\frac14 \beta_n) \bar B_{4n-2},
\;\; B_{2n+1}=\frac12\bar B_{4n}+\bar B_{4n+1}.
$$

\pf
We only give the proof for $B_{2i}$, the others being similar. Let
$$
B_{2i}=\mu_{4i-2}\bar B_{4i-2}+\mu_{4i-1}\bar B_{4i-1}+\mu_{4i}\bar
B_{4i}+\mu_{4i+1}\bar
B_{4i+1}.
$$
By application of the CCA, starting from the B-coefficients
$[0,0,1,\omega_i]$ of $B_{2i}$
on  the interval $I_i$, we deduce the B-coefficients of $B_{2i}$,
respectively on the
subintervals $I'_i$ and $I''_i$:
$$
[0,0, \frac12+\frac14\beta_i, \frac12] \;\; and \;\;
[\frac12, \frac12-\frac14\beta_i, \frac12+\frac12 \omega_i, \omega_i].
$$
Similarly, from the B-coefficients $[\omega_i,0,0,0]$ of $B_{2i}$ on
$I_{i+1}$, we deduce the
B-coefficients of $B_{2i}$ on the subintervals $I'_{i+1}$ and $I''_{i+1}$:
$$
[\omega_i,\frac12\omega_i,0,0] \;\; and \;\; [0,0,0,0]
$$
On the other hand, the B-coefficients of the finer B-splines are respectively
\ms
\noin
1) on the subintervals $I'_i$ and $I''_i$
$$
for \;\; \bar B_{4i-2}: \;\; [0,0,1,\frac12], \;\;[\frac12,0,0,0]
$$
$$
for \;\;\bar B_{4i-1} :\;\; [0,0,0,\frac12], \;\;[\frac12,1,0,0]
$$
2) on the subintervals $I''_i$ and $I'_{i+1}$:
$$
for \;\;\bar B_{4i}: \;\; [0,0,1,\omega_i], \;\;[\omega_i,0,0,0]
$$
$$
for \;\;\bar B_{4i+1}: \;\; [0,0,0,1-\omega_i], \;\;[1-\omega_i,1,0,0]
$$
Therefore, the B-coefficients of $B_{2i}$ as linear combination of the four
finer B-splines
on the three subintervals $I'_i$, $I''_i$ and $I'_{i+1}$ are respectively
equal to
$$
[0,0,\mu_{4i-2},\frac12(\mu_{4i-2}+\mu_{4i-1})],\;\;
[\frac12(\mu_{4i-2}+\mu_{4i-1}),\mu_{4i-1},\mu_{4i},\omega_i\mu_{4i}+(1-\omega_i)\mu_{4i+1}]
$$
$$
[\omega_i\mu_{4i}+(1-\omega_i)\mu_{4i+1},\mu_{4i+1},0,0].
$$
By identifying these B-coefficients with those of $B_{2i}$, one obtains:
$$
\mu_{4i-2}=\frac12+\frac14 \beta_i, \;\mu_{4i-1}=\frac12-\frac14 \beta_i,
\;\mu_{4i}=\frac12 (1+\omega_i), \; \mu_{4i+1}=\frac12 \omega_i
$$

\eop

\noin
An immediate consequence of the previous theorem is the following {\sl
global corner-cutting algorithm}:

\proclaim Theorem 7.
Given the two expansions of $S\in GS_2(I, X,\beta)$ with respect to the
coarse and fine
B-splines bases:
$$
S=\sum_{k=0}^{2n+1} \gamma_k B_k=\sum_{l=0}^{4n+1} \delta_l \bar B_l,
$$
then, the new B-coefficients have the following expressions in terms of the
former B-coefficients:
$$
\delta_{4i-2}=(\frac12-\frac14\beta_i)\gamma_{2i-1}+(\frac12+\frac14\beta_i)
\gamma_{2i},
$$
$$
\delta_{4i-1}=(\frac12+\frac14\beta_i)\gamma_{2i-1}+(\frac12-\frac14\beta_i)
\gamma_{2i},
$$
$$
\delta_{4i}=\frac12(1+\omega_i)\gamma_{2i}+\frac12(1-\omega_i)\gamma_{2i+1},
$$
$$
\delta_{4i+1}=\frac12\omega_i\gamma_{2i}+\frac12(2-\omega_i)\gamma_{2i+1}.
$$

\pf
The proof simply consists in comparing the coefficients in the two
expressions after
substituting in the first expression the coarser B-splines $B_k$ by their
expansions as linear
combinations of the finer B-splines $\bar B_l$, given in theorem 3.
\eop

\proclaim Theorem 8.
 The sequence of SCPs associated with a given GQS, obtained by successive
applications of the global CCA,
 converges uniformly to the GQS.

\pf
Let $S(x)=\sum_{k=0}^{2n+1}\gamma_k B_k$ be the equation of a GQS and let
$P_0=\sum_{k=0}^{2n+1} \gamma_k \phi_k$ be the equation of its initial SCP.
We denote by $\phi_{2i}$ the hat function with support
$[\eta_{i-1},\eta_i]$ satisfying $\phi_{2i}(\xi_i)=1$, and by $\phi_{2i+1}$
the hat function with support
$[\xi_{i},\xi_{i+1}]$ satisfying $\phi_{2i}(\eta_i)=1$.
\ms
\noin
For $x\in [\eta_{i-1},\xi_i]$, (resp. $x\in [\xi_i,\eta_i]$), we have
$$
\vert S(x)-P_0(x)\vert \le \max \{\vert S(x)-\gamma_{2i-1}\vert, \vert
S(x)-\gamma_{2i}\vert\}
$$
$$
(resp. \vert S(x)-P_0(x)\vert \le \max \{\vert S(x)-\gamma_{2i}\vert, \vert
S(x)-\gamma_{2i+1}\vert\}).
$$
Moreover, for $x\in [\eta_{i-1},\xi_i]$, we have
$$
\vert S(x)-\gamma_{2i}\vert \le \sum_{k=2i-2}^{2i+1}\vert
\gamma_k-\gamma_{2i}\vert B_k(x)\le \max \{\vert \gamma_k-\gamma_{2i}\vert,
2i-2\le k\le 2i+1\}
$$
As we observe that
$$
\vert \gamma_k-\gamma_{2i}\vert \le 2\max \{\vert
\gamma_{k+1}-\gamma_k\vert, 2i-2\le k\le 2i+1\},
$$
we are led to define
$$
\Delta_0=max \{\vert \gamma_{k+1}-\gamma_k\vert, 0\le k\le 2n\},
$$
and we finally obtain
$$
\Vert S-P_0\Vert_{\infty}\le 2 \Delta_0.
$$
Now, we do the same for the next SCP $P_1$ deduced from $P_0$ by one
application of the global CCA.
Here, we define
$$
\Delta_1=\max \{\vert \delta_{l+1}-\delta_l\vert, 0\le l\le 4n+1\}.
$$
We have successively the following majorations
$$
\vert \delta_{4i-1}-\delta_{4i-2}\vert=\vert \frac12
\beta_i\;(\gamma_{2i}-\gamma_{2i-1})\vert \le \frac12\vert
\gamma_{2i}-\gamma_{2i-1}\vert
$$
$$
\vert \delta_{4I+1}-\delta_{4i}\vert=\frac12 \vert
\gamma_{2i+1}-\gamma_{2i}\vert
$$
$$
\vert \delta_{4i}-\delta_{4i-1}\vert=\vert
\frac12(1-\omega_i)(\gamma_{2i+1}-\gamma_{2i})+(\frac12+\frac14\beta_i)(\gamma_{2i}-\gamma_{2i-1}))\vert
$$
$$
\le \frac12\max\{\vert \gamma_{2i+1}-\gamma_{2i}\vert, \vert
\gamma_{2i}-\gamma_{2i-1}\vert\}
$$
Therefore we obtain
$$
\Delta_1\le \frac12 \Delta_0.
$$
Denoting by $\Delta_m$ the maximum distance between two consecutive
vertices of the SCP $P_m$
obtained after $m$ applications of the global CCA, we get
$$
\Delta_m\le \frac{1}{2^m}\Delta_0,
$$
which proves the uniform convergence to the GQS $g$ of the sequence $(P_n)$
of its SCPs.
\eop
*******************

\sect{6. Quasi-interpolant and Lagrange interpolant}

\noin
In this section, we define two kinds of approximation operators: the first
is a quasi-interpolant
with good shape-preserving properties while the second is a Lagrange
interpolant which is uniformly
bounded independently of the partition.
Let us begin with the {\sl quasi-interpolant} $Q$ defined by
$$
Qf=\sum_{i=0}^{n+1} [f(\xi_i)B_{2i}+f(\eta_i) B_{2i+1}]
$$

\proclaim Theorem 9.
$Q$ is exact on $\PP_1$, $\Vert Q \Vert_{\infty}=1$,
and $Q$ preserves the monotonicity and the convexity of $f$.

\pf
Since $e_0=\sum_{i=0}^{n+1} [B_{2i}+B_{2i+1}]$ and $e_1=\sum_{i=0}^{n+1}
[\xi_i B_{2i}+\eta_i B_{2i+1}]$,
we obtain immediately the property that $Q$ is exact on $\PP_1$. Moreover,
$Qe_0=e_0$ and
$\vert Qf \vert_{\infty}\le \vert f\vert_{\infty}$ imply $\Vert Q
\Vert_{\infty}=1$.

\noin
If $f$ is increasing (resp. convex), the global B-polygon of $Qf$ is also
increasing (resp. convex),
and the result follows by a direct application of theorem 2.
\eop

\noin
Now, let us study a {\sl Lagrange interpolation operator} $L$ in the space
$GS_2(\beta)$.
Let $m'_i$ (resp. $m''_i$) be the midpoint of $I'_i=[x_{i-1},m_i]$ (resp.
of $I''_i=[m_i,x_i]$.
The following theorem is an extension to GQS of a result previously given by
Kammerer et al. [4] and by Marsden [5] for ordinary quadratic splines.

\proclaim Theorem 10.
(i) Given a function $f$ defined on $I=[a,b]$, there exists a unique
generalized quadratic spline
$Lf\in GS_2(I,X,\beta)$  satisfying the following interpolation properties:
$$
Lf(a)=f(a), \;\; Lf(b)=f(b)\;\; and \;\;for \;\; 1\le i\le n,
$$
$$
Lf(m'_i)=f(m'_i)\;\; and \;\; Lf(m''_i)=f(m''_i).
$$
(ii) The B-coefficients of $Lf=\sum_{k=0}^{2n+1} \gamma_k B_k$ are
solutions of the tridiagonal system of $2n$
linear equations ($1\le i\le n$):
$$
\omega_{i-1}\gamma_{2i-2}+(3-\omega_{i-1}-\frac12\beta_i)
\gamma_{2i-1}+(1+\frac12\beta_i) \gamma_{2i}=4f(m'_i)
$$
$$
(1+\frac12\beta_i) \gamma_{2i-1}+(2+\omega_i-\frac12\beta_i)
\gamma_{2i}+(1-\omega_i) \gamma_{2i+1}=4f(m''_i)
$$
with $\gamma_0=f(a)$, $\gamma_{2n+1}=f(b)$.

\pf
On the interval $I_i$, we are able to compute the values of
$Lf=\sum_{k=2i-2}^{2i+1} \gamma_k B_k$ at the
points $m'_i$ and $m''_i$ by applying the CCA to B-splines. Without going
into details, it is straightforward to
obtain the coefficients of the two equations given in the theorem. The
tridiagonal matrix of the linear system
is strictly diagonally dominant since the two corresponding inequalities
$$
3-\omega_{i-1}-\frac12\beta_i>\omega_{i-1}+1+\frac12\beta_i \;\; and \;\;
2+\omega_i-\frac12\beta_i>1+\frac12\beta_i+1-\omega_i
$$
are respectively equivalent to the following
$$
2(1-\omega_i)-\beta_i>0 \;\; and \;\; 2\omega_i-\beta_i>0
$$
which are obviously satisfied since $\beta_i<0$ and $0<\omega_i<1$.
\eop

\proclaim Theorem 11.
 $\Vert L\Vert_{\infty}$ is uniformly bounded for all partitions of $I$.
More specifically,
setting $\bar\beta=max\{\beta_i,1\le i\le n\}$, we obtain
$$
\Vert L \Vert_{\infty} \le \frac{4(3\bar\beta-1)}{\bar\beta(5-3\bar\beta)}.
$$

\pf
For this purpose, we start from the expression of $Lf=\sum_{r=0}^{3}
a_{i,r} b_r$ in the local B-spline
basis of each subinterval $I_i$. Then, by
the CCA, we compute the values of the
$b_r$ at the two points $m'_i$ and $m''_i$ and we get the two equations:
$$
4Lf(m'_i)=a_{i,0}+(2-\frac12\beta_i)a_{i,1}+(1+\frac12\beta_i)a_{i,2}=4f(m'_i)
$$
$$
4Lf(m''_i)=(1+\frac12\beta_i)a_{i,1}+(2-\frac12\beta_i)a_{i,2}+a_{i,3}=4f(m'_i)
$$
We have to add to these equations the $C^1$ continuity conditions of $Lf$
at interior points $x_i$,
which can be written:
$$
a_{i-1,3}=a_{i,0}\;\; and \;\;
\lambda_{i-1}(a_{i-1,3}-a{i-1,2})=\lambda_i(a_{i,1}-a_{i,0})
$$
where
$\lambda_i=2\left(\frac{\beta_i-1}{\beta_i}\right)=\frac{1}{\theta_i}$ for
all $1\le i\le n$.

\noin
>From these conditions, we deduce the expressions
$$
a_{i-1,3}=a_{i,0}=\frac{\lambda_{i-1}a_{i-1,2}+\lambda_i
a_{i,1}}{\lambda_{i-1}+\lambda_i}
$$
that we substitute in the two equations above. Taking as unknowns
$a_{2i-1}=a_{i,1}$ and
$a_{2i}=a_{i,2}$, we then obtain the following tridiagonal system of equations
$$
\frac{\lambda_{i-1}}{\lambda_{i-1}+\lambda_i}a_{2i-2}+\left[\frac{\lambda_i}
{\lambda_{i-1}+\lambda_i}+2-\frac12\beta_i\right]a_{2i-1}+(1+\frac12\beta_i)
a_{2i}=4f(m'_i)
$$
$$
(1+\frac12\beta_i)a_{2i-1}+\left[\frac{\lambda_i}{\lambda_i+\lambda_{i+1}}+2
-\frac12\beta_i\right]a_{2i}+\frac{\lambda_{i+1}}{\lambda_i+\lambda_{i+1}}a_
{2i+1}=4f(m''_i)
$$
The tridiagonal matrix is strictly diagonally dominant since the two conditions
$$
\frac{\lambda_i}{\lambda_{i-1}+\lambda_i}+2-\frac12\beta_i>\frac{\lambda_{i-
1}}{\lambda_{i-1}+\lambda_i}+1+\frac12\beta_i,
$$
$$
\frac{\lambda_i}{\lambda_i+\lambda_{i+1}}+2-\frac12\beta_i>1+\frac12\beta_i+
\frac{\lambda_{i+1}}{\lambda_i+\lambda_{i+1}},
$$
can also be written
$$
-\beta_i>\frac{-2\lambda_i}{\lambda_{i-1}+\lambda_i} \;\; and \;\;
-\beta_i>\frac{-2\lambda_i}{\lambda_{i+1}+\lambda_i},
$$
and they are obviously satisfied since $\beta_i<0$ and $\lambda_i>0$ for
all $i$.

\noin
Let $\vert a \vert_{\infty}=max\{\vert a_k \vert, 0\le k\le 2n+1\}$.
>From the above tridiagonal system, we deduce respectively  the following
inequalities
$$
\left[\frac{\lambda_i}{\lambda_{i-1}+\lambda_i}+2-\frac12\beta_i\right]\vert
a_{2i-1} \vert \le 4\vert f(m'_i) \vert +
\frac{\lambda_{i-1}}{\lambda_{i-1}+\lambda_i}\vert a_{2i-2} \vert
+(1+\frac12\beta_i)\vert a_{2i}
\vert,
$$
$$
\left[\frac{\lambda_i}{\lambda_i+\lambda_{i+1}}+2-\frac12\beta_i\right]\vert
a_{2i} \vert \le 4\vert f(m''_i) \vert +
(1+\frac12\beta_i)\vert a_{2i-1}
\vert+\frac{\lambda_{i+1}}{\lambda_i+\lambda_{i+1}}\vert a_{2i+1}
\vert,
$$
and we obtain for all $i$:
$$
\left[\frac{2\lambda_i}{\lambda_{i-1}+\lambda_i}-\beta_i\right]\vert a
\vert_{\infty}\le 4\Vert f \Vert_{\infty}
\;\; and \;\;
\left[\frac{2\lambda_i}{\lambda_{i+1}+\lambda_i}-\beta_i\right]\vert a
\vert_{\infty}\le 4\Vert f
\Vert_{\infty}.
$$
Defining $\bar\lambda=max\{\lambda_i, 1\le i\le n\}$ and
$\bar\beta=max\{\beta_i,1\le i\le n\}$, then
we have $4\le \lambda_i\le\bar\lambda=2\frac{\bar\beta-1}{\bar\beta}$, hence
$$
\frac{2\lambda_i}{\lambda_{i-1}+\lambda_i}-\beta_i\ge
\frac{8}{\lambda_{i-1}+4}-\bar\beta\ge\frac{8}{\bar\lambda+4}-\bar\beta=\frac{\bar\beta(5-3\bar\beta)}{(3\bar\beta-1)}
$$
Therefore, we finally obtain
$$
\vert a \vert_{\infty}\le \frac{4(3\bar\beta-1)}{\bar\beta(5-3\bar\beta)}
\Vert f \Vert_{\infty}.
$$
As the graph of $Lf$ lies in the convex hull of its local control polygons,
it is not difficult to see that
we have also
$$
\Vert Lf \Vert_{\infty}\le \frac{4(3\bar\beta-1)}{\bar\beta(5-3\bar\beta)}
\Vert f \Vert_{\infty}.
$$
In other words, we obtain the following upper bound which is independent of the given partition, but only depends
on the sequence $\beta$:
$$
\Vert L \Vert_{\infty} \le \frac{4(3\bar\beta-1)}{\bar\beta(5-3\bar\beta)}.
$$
\eop

\noin
{\bf Remarks}: 1. In the case of classical quadratic splines, we have
$\bar\beta=-1$, whence $\Vert L
\Vert_{\infty} \le 2$, which is the sharp upper bound given in [4] and [5].
This suggests that the
upper bound given in theorem 11 may also be sharp.

\ms
\noin
2. As for error bounds, we know, from a classical result in approximation
theory (see [3], theorem
), that for $P=Q$ or $L$
$$
\Vert f-Pf\Vert_{\infty}\le (1+\Vert P\Vert_{\infty})d_{\infty}(f,GS_2(\beta)).
$$
For  $f\in C^2(I)$, we already know (see [3]) that
$d_{\infty}(f,GS_2(\beta))=O(h^2)$ ,
with $h=\max h_i$. Therefore, we see that both operators have an
approximation order equal to
$2$.

*********************************************
\References

\ref
P. Costantini, On monotone and convex spline interpolation. Math. of
Comput. {\bf 46}, No 173
(1986), 203-214.

\ref
P. Costantini, Co-monotone interpolating splines of arbitrary degree. A
local approach.
SIAM J. Sci. Stat. Comput.  {\bf 8}, No 6 (1987), 1026-1034.

\ref
R.A. DeVore, G.G. Lorentz, {\sl Constructive Approximation},
Springer-Verlag, Berlin, 1993.

\ref
Kammerer, Reddien, R.S. Varga, Quadratic interpolatory splines, Numer.
Math. {\bf 22} (1974),
241-259.

\ref
M. Marsden, Operator norm bounds and error bounds for quadratic spline
interpolation.
In {\sl Approximation Theory}, Banach Center Publications, vol. {\bf 4}
(1979), 159-175.

\ref
J.-L. Merrien, A family of Hermite interpolants by bisection
algorithms. \NA \ {\bf 2} (1992), 187-200.

\ref
J.L. Merrien, P. Sablonni\`ere, Monotone and convex $C^1$ Hermite interpolants
generated by a subdivision algorithm, \CA \ {\bf 19} (2003), 279-298.

\ref
J.L. Merrien, P. Sablonni\`ere, Monotone and convex $C^1$ Hermite interpolants
generated by an adaptive subdivision scheme, \CRASP \ t. 333, S\'erie I, p.
493-497,
2001.

\ref
C.A. Micchelli, {\sl Mathematical Methods in CAGD}. SIAM, Philadelphia, 1995.

\ref
P. Sablonni\`ere, Bernstein type bases and corner-cutting algorithms for
$C^1$ Merrien curves.
In Proceedings of the conference {\sl Multivariate approximation and
interpolation},
Almu\~{n}ecar, Spain (2001)). Adv. in Comput. Math. {\bf 20} (2004), 229-246.

\ref
I.J. Schoenberg, {\sl Selected papers}, Volumes 1 and 2, edited by C. de Boor.
Birkh\"auser-Verlag, Boston, 1988.

\ref
L.L. Schumaker, {\sl Spline functions: basic theory}. John Wiley and Sons,
New-York, 1981.


\Address

 P. Sablonni\`ere, INSA de Rennes,
20 Avenue des Buttes de Co\"esmes,

 CS 14315, 35043 RENNES CEDEX, France
{\tt psablonn@insa-rennes.fr}

\bye